# A Note on the Forward-Douglas–Rachford Splitting for Monotone Inclusion and Convex Optimization

**Hugo Raguet**
Aix-Marseille University, *I2M*
hugo.raguet@gmail.com

### Abstract

We shed light on the structure of the "three-operator" version of the forward-Douglas–Rachford splitting algorithm for finding a zero of a sum of maximally monotone operators $A + B + C$, where $B$ is cocoercive, involving only the computation of $B$ and of the resolvent of $A$ and of $C$, separately. We show that it is a straightforward extension of a fixed-point algorithm proposed by us as a generalization of the forward-backward splitting algorithm, initially designed for finding a zero of a sum of an arbitrary number of maximally monotone operators $\sum_{i=1}^{n} A_i + B$, where $B$ is cocoercive, involving only the computation of $B$ and of the resolvent of each $A_i$ separately. We argue that, the former is the "true" forward-Douglas–Rachford splitting algorithm, in contrast to the initial use of this designation in the literature. Then, we highlight the extension to an arbitrary number of maximally monotone operators in the splitting, $\sum_{i=1}^{n} A_i + B + C$, in a formulation admitting preconditioning operators. We finally demonstate experimentally its interest in the context of nonsmooth convex optimization.

## 1 Introduction and Motivation

We introduced some time ago a *generalization of the forward-backward* splitting algorithm (Raguet et al., 2013) for solving, over a real Hilbert space $\mathcal{H}$, monotone inclusion problems of the form

$$\text{find } x \in \text{zer}\left( \sum_{i=1}^{n} A_i + B \right), \tag{P1}$$

by making use only of the *resolvent* of each set-valued operator $A_i \colon \mathcal{H} \to 2^{\mathcal{H}}$, supposed maximally monotone, and of the explicit application of $B \colon \mathcal{H} \to \mathcal{H}$, supposed *cocoercive*. This task is especially interesting when identifying operators to subdifferentials; under suitable conditions, it is equivalent to

$$\text{find } x \in \arg\min \sum_{i=1}^{n} g_i + f, \tag{P2}$$

by making use only of the *proximity operator* of each proper, convex, and lower semicontinuous functionals $g_i \colon \mathcal{H} \to \left]-\infty, +\infty\right]$, and of the gradient of $f \colon \mathcal{H} \to \mathbb{R}$, supposed differentiable with a Lipschitz-continuous gradient.

To this end, we recast the problem on a convenient augmented space $\boldsymbol{\mathcal{H}} \stackrel{\text{def}}{=} \mathcal{H}^n$ as finding a zero of the sum $\boldsymbol{A} + \boldsymbol{B} + \boldsymbol{N}_\mathcal{S}$, where the maximally monotone operators $\boldsymbol{A}$ and $\boldsymbol{B}$ depend on each $A_i$ and on $B$, respectively, and $\boldsymbol{N}_\mathcal{S}$ is the *normal cone* of a suitable subspace $\mathcal{S}$ (see § 2.2). Then, we solve the problem using a fixed-point algorithm involving only the above mentioned resolvents and operator together with the resolvent of $\boldsymbol{N}_\mathcal{S}$, which is nothing but the *orthogonal projector* onto $\mathcal{S}$. The resulting method can be viewed as a combination of the well-known *forward-backward* and *Douglas–Rachford* splitting algorithms, recalled in § 2.1.

Actually, the exact same fixed-point algorithm can be used for finding a zero of a sum $A + B + N_\mathcal{V}$, where now $N_\mathcal{V}$ is the normal cone of any closed subspace in $\mathcal{H}$ (we revert to nonbold notations, since the space on which one solves the inclusion is general); the above mentioned orthogonal projection is now onto $\mathcal{V}$. This has been noted by Briceño-Arias (2015); although the practical use of this improvement is limited, it gives some theoretical insight. Due to the structure of the algorithm, it is coined *forward-Douglas–Rachford*.

However, in our opinion the work is not done yet for such a designation. As advocated below, the complete forward-Douglas–Rachford should be able to find a zero of a sum $A + B + C$, where now $C$ is any maximally monotone operator. As simplistic as it might sound, this is done by simply replacing the orthogonal projection onto $\mathcal{V}$ by the resolvent of $C$. A subtle difference is due to the fact that this resolvent is not linear anymore, requiring some more work for deriving the full convergence analysis. This has been finally done by Davis and Yin (2017), who strangely fail to mention our work. Note that these authors also study rates of convergence, and possible acceleration schemes. Let us also mention the stochastic version developed by Cevher et al. (2016).

In this paper, we want to precisely describe the links between the above methods; this is done in § 2.

Then, since our interest lies mostly in large-scale convex optimization problems involving many terms, we explicitly derive the forward-Douglas–Rachford splitting algorithm for solving the monotone inclusion problem

$$\text{find } x \in \text{zer}\left( \sum_{i=1}^n A_i + B + C \right); \tag{P3}$$

as far as we know, this has not been written yet. Problem P3 is the same as problem P1 with an additional maximally monotone operator $C$. Observe that, although the operator $C$ satisfies the same assumptions as each $A_i$, we particularize its role in the splitting of problem P3. In the resulting iterative method, this will translate to two desirable properties which we describe below. Note that we are not aware of a situation where these properties are crucial, though, and in practice any problem that can be solved by the forward-Douglas–Rachford splitting algorithm can also be solved by the generalized forward-backward; but in the following cases, the former is somewhat more elegant.

The first property is that the resolvent of $C$ is directly applied to the iterates. One consequence is that the iterates always lies within the *domain* of $C$ along iterations. Hence, the algorithm can be applied when $B$ is cocoercive only within this domain, without further care. Moreover, this better handles a convex hard constraint, i.e. when $C$ is the normal cone of a closed convex set, since such constraint will be satisfied at each iteration. In contrast, the generalized forward-backward ensures convex hard constraints only at convergence; they might be violated during iterations. More generally, this is useful for instance when $C$ is the subdifferential of the $\ell_1$-norm in some basis, inducing *sparsity* (many coefficients in that basis are zero) over the solution set of a convex optimization problem. The resolvent of $C$ is the so-called *soft-thresholding*, squeezing nonsignificant coefficients to zero. Again, the generalized forward-backward ensures sparsity only at convergence, and coefficients which are zero at a solution may be nonzero (with decreasing amplitude) during iterations; while the iterates in the forward-Douglas–Rachford might have the right support after a finite number of iterations. In turn, this could be useful for example for risk estimation, (see Zou et al., 2007) ; see also our numerical experiments § 4.

The second property is that no additional auxiliary variable needs to be stored in memory for taking into account the operator $C$ in the splitting. More precisely, in principle $n + 1$ auxiliary variables, each of the dimension of the problem, must be stored for applying the generalized forward-backward splitting algorithm to problem P3, against only $n$ auxiliary variables for applying forward-Douglas–Rachford as shown in § 3. However, practical implementations usually has other memory needs, so that the memory gained this way rarely exceeds one fifth of the overall memory needed, even when $n$ is as low as one.

Let us precise that there exists now many methods able to solve problem P3 while taking advantage of both the splitting and the cocoercivity of $B$, in the sense that the required operator inversions involve only each $A_i$ and $C$ individually, and not $B$. Most notable examples in the literature either follow the primal-dual approach proposed independently by Condat (2013) and Vũ (2013), or the one of Combettes and Pesquet (2012), or our method (Raguet et al., 2013). Note that the mentioned primal-dual methods can deal with a larger class of problems. Moreover with the former, it is easy to particularize the role of one of the maximally monotone operators in the splitting, with the desirable algorithmic properties explained above. Nonetheless, explicit primal-dual algorithms are not necessarily best suited for the considered class of problems for two reasons. First, they allow only for restricted range of parameters (namely the explicit step size, denoted hereafter $\gamma$ and depending on the cocoercivity constant of $B$, see Raguet, 2014, III.2.3), and second, our setting can be more appropriate for preconditioning purpose (see our work on this topic Raguet and Landrieu, 2015).

For those reasons, we specify in § 3 the algorithm for both monotone inclusion and convex optimization problems structured as problem P3, with the possibility of using preconditioners. Finally, § 4 is devoted to numerical illustration.

## 2   The "True" Forward-Douglas–Rachford Operator

As we will see, both the exposition of our method and the discussion of its relationship to others, are facilitated within the unifying framework of Combettes (2004). We also refer the less familiar readers to this article for most specific notations used in the following.

### 2.1   Forward-Backward and Douglas–Rachford Algorithms viewed as Compositions of Averaged Operators

The forward-backward splitting algorithm is a well-studied method for solving monotone inclusion problem P1 when $n$ is restricted to one, that is

$$\text{find } x \in \text{zer}(A + B) \, ,$$

by making use only of the resolvent of $A$ and of the explicit application of $B$.

Combettes (2004, Section 6) gives both synthetic description and analysis, using the fact that the set of *zeros* of the sum of $A + B$ is equal to the set of *fixed points* of the operator

$$T_{\text{FB}} \stackrel{\text{def}}{=} J_{\gamma A}(\text{Id} - \gamma B) \, , \tag{2.1}$$

where $\gamma \in \mathbb{R}_{++}$ is a positive scalar parameter, Id is the identity over $\mathcal{H}$, and $J_{\gamma A} \stackrel{\text{def}}{=} (\text{Id} + \gamma A)^{-1}$ is the resolvent of $\gamma A$. Then, the algorithm is the repeated application of $T_{\text{FB}}$; since $J_{\gamma A}$ is *firmly nonexpansive* and, for well chosen $\gamma$, $(\text{Id} - \gamma B)$ is also $\alpha$-*averaged* for some $\alpha$, so is their composition. Weak convergence towards a fixed point ensues thanks to results given in the cited paper, allowing for relaxations, varying step size $\gamma$ and summable error terms.

Alternatively, the Douglas–Rachford splitting algorithm is an equally well-studied method, for solving monotone inclusion problem P3 when $n$ is restricted to one and $B$ is restricted to zero, that is

$$\text{find } x \in \text{zer}(A + C),$$

involving only the resolvents of $A$ and $C$, separately.

Again, Combettes (2004, Section 5) gives both synthetic description and analysis, using now the fact that the *preimage* of the set of zeros of the sum of $A + C$ by the single-valued operator $J_{\gamma C}$ is equal to the set of fixed points of the operator $R_{\gamma A} R_{\gamma C}$, where $R_{\gamma A} \stackrel{\text{def}}{=} 2 J_{\gamma A} - \text{Id}$ is the *reflection* operator associated to $\gamma A$. Here, it is convenient to observe that $R_{\gamma A} R_{\gamma C}$ shares its fixed points with the operator

$$T_{\text{DR}} \stackrel{\text{def}}{=} \tfrac{1}{2}(R_{\gamma A} R_{\gamma C} + \text{Id}) = J_{\gamma A}(2 J_{\gamma C} - \text{Id}) + (\text{Id} - J_{\gamma C}). \tag{2.2}$$

Then, the algorithm is usually described as the repeated application of the latter operator, which is firmly nonexpansive for any $\gamma \in \mathbb{R}_{++}$. Weak convergence towards a fixed point ensues thanks to results given in the cited paper, allowing for relaxations and essentially summable error terms. Finally, a zero is found by applying $J_{\gamma C}$ to such fixed point; in contrast to the above, varying step size $\gamma$ along iterations is not possible in general, since the set of fixed points would vary as well.

Let us mention that a proof of the firm nonexpansiveness of the operator defined above goes back to the work of Lions and Mercier (1979), and introduction of relaxation and summable error terms is due to Eckstein and Bertsekas (1992).

## 2.2   Generalized Forward-Backward and Forward-Douglas–Rachford Algorithms

In terms of convex optimization, the Douglas–Rachford algorithm can be applied to minimize a sum of two convex functionals using their proximity operators. If one needs to split the objective into more than two functionals, a useful trick is that finding $x \in \mathcal{H}$ minimizing $\sum_{i=1}^{n} g_i(x)$ is mathematically equivalent to finding $x_1, \ldots, x_n \in \mathcal{H}$ minimizing $\sum_{i=1}^{n} g_i(x_i)$ under the constraint $x_1 = \cdots = x_n$. That is, introducing the augmented space $\boldsymbol{\mathcal{H}} \stackrel{\text{def}}{=} \mathcal{H}^n$, the closed subspace $\boldsymbol{\mathcal{S}} \stackrel{\text{def}}{=} \{\boldsymbol{x} = (x_1, \ldots, x_n) \in \boldsymbol{\mathcal{H}} \mid x_1 = \cdots = x_n\}$ and its convex indicator $\iota_{\boldsymbol{\mathcal{S}}} \colon \boldsymbol{\mathcal{H}} \to ]-\infty, +\infty] \colon \boldsymbol{x} \mapsto 0$ if $\boldsymbol{x} \in \boldsymbol{\mathcal{S}}$, $+\infty$ otherwise, the problem is to find $\boldsymbol{x} \in \arg\min \boldsymbol{g} + \iota_{\boldsymbol{\mathcal{S}}}$, where the functional $\boldsymbol{g} \colon \boldsymbol{\mathcal{H}} \to ]-\infty, +\infty] \colon \boldsymbol{x} \mapsto \sum_{i=1}^{n} g_i(x_i)$ is decoupled along the splitting. Its proximity operator can be deduced from the proximity operators of each $g_i$ separately, and the proximity operator of $\iota_{\boldsymbol{\mathcal{S}}}$ is nothing but the orthogonal projector onto $\boldsymbol{\mathcal{S}}$, essentially averaging the components.

This approach translates easily for solving monotone inclusions involving the sum of an arbitrary number of maximally monotone operators. Introducing the product operator $\boldsymbol{A} \stackrel{\text{def}}{=} \times_{i=1}^{n} A_i \colon \boldsymbol{\mathcal{H}} \to 2^{\boldsymbol{\mathcal{H}}}$, and the normal cone $\boldsymbol{N}_{\boldsymbol{\mathcal{S}}} \colon \boldsymbol{\mathcal{H}} \to 2^{\boldsymbol{\mathcal{H}}} \colon \boldsymbol{x} \mapsto \boldsymbol{\mathcal{S}}^{\perp}$ if $\boldsymbol{x} \in \boldsymbol{\mathcal{S}}$, $\emptyset$ otherwise, one can show that finding $x \in \text{zer}(\sum_{i=1}^{n} A_i x)$ is equivalent to finding $\boldsymbol{x} \in \text{zer}(\boldsymbol{A} + \boldsymbol{N}_{\boldsymbol{\mathcal{S}}})$. Again, the Douglas–Rachford algorithm can be applied with full splitting since the resolvent of $\boldsymbol{A}$ can be deduced from the resolvents of each $A_i$ separately, and the resolvent of $\boldsymbol{N}_{\boldsymbol{\mathcal{S}}}$ is again the orthogonal projector onto $\boldsymbol{\mathcal{S}}$. The resulting algorithm (with the step size $\gamma$ restricted to unity) coincide with the *method of* *partial inverse* of Spingarn (1983), who fails to notice the connection with the Douglas–Rachford algorithm.

Oddly enough, it is only thirty years later that one could combine such a splitting into an arbitrary number of maximally monotone operator, together with an additional operator handled through its explicit application rather than its resolvent; that is, solve our problem P1. In a nutshell, we introduce also the operator $\boldsymbol{B}\colon\mathcal{H}\to\mathcal{H}\colon \boldsymbol{x}\mapsto(Bx_1,\ldots,Bx_n)$, cast the problem into finding $\boldsymbol{x}\in\mathrm{zer}(n\boldsymbol{A}+\boldsymbol{B}+\boldsymbol{N}_{\mathcal{S}})$, the latter set being the orthogonal projection onto $\mathcal{S}$ of the set of fixed points of the composition

$$\boldsymbol{T}_{\mathrm{GFB}} \stackrel{\mathrm{def}}{=} \tfrac{1}{2}\bigl(\boldsymbol{R}_{n\gamma \boldsymbol{A}}\boldsymbol{R}_{\mathcal{S}}+\mathrm{Id}\bigr)\bigl(\mathrm{Id}-\gamma \boldsymbol{B}\boldsymbol{P}_{\mathcal{S}}\bigr) = \boldsymbol{J}_{n\gamma \boldsymbol{A}}\bigl(2\boldsymbol{P}_{\mathcal{S}}-\mathrm{Id}-\gamma \boldsymbol{B}\boldsymbol{P}_{\mathcal{S}}\bigr) + \bigl(\mathrm{Id}-\boldsymbol{P}_{\mathcal{S}}\bigr)\ . \qquad(2.3)$$

This operator is obviously a combination of both the forward-backward equation 2.1 and the Douglas–Rachford equation 2.2 operators in the augmented space setting, where we particularize the orthogonal projector $\boldsymbol{P}_{\mathcal{S}} \stackrel{\mathrm{def}}{=} \boldsymbol{J}_{\boldsymbol{N}_{\mathcal{S}}}$ and the reflection operator $\boldsymbol{R}_{\mathcal{S}} \stackrel{\mathrm{def}}{=} 2\boldsymbol{P}_{\mathcal{S}}-\mathrm{Id}$. We use the fact that, provided that $B$ is cocoercive in $\mathcal{H}$, then so is the composition $\boldsymbol{P}_{\mathcal{S}}\boldsymbol{B}\boldsymbol{P}_{\mathcal{S}} = \boldsymbol{B}\boldsymbol{P}_{\mathcal{S}}$ in $\mathcal{H}$. Thanks to the direct composition appearing in the middle term in equalities equation 2.3, it is straightforward to show that it is $\alpha$-averaged, and then to prove convergence towards a fixed point from repeated applications with possible relaxations and summable errors.

Since our method resembles a forward-backward algorithm on an augmented space, we coin it a *generalized forward-backward algorithm*, and exemplify in our article (Raguet et al., 2013) its usefulness for large-scale convex optimization problems. Shortly after, Briceño-Arias (2015) notices that the vector space $\mathcal{S}$ can actually be replaced by any closed vector space $\mathcal{V}$ without modifying the whole analysis, and writes the resulting method for finding (we revert again to nonbold notations for generality) $x\in\mathrm{zer}(A+B+N_{\mathcal{V}})$. The resulting algorithm is slightly different because the identity $P_{\mathcal{V}}BP_{\mathcal{V}} = BP_{\mathcal{V}}$ is not true in general; for some applications, it could trade some memory savings against increased computational load, without significant difference. The author also revisits the approach of the method of partial inverse in this setting, without practical applications known to us.

Finally, it is now easy to see a more interesting extension of our generalized forward-backward algorithm, consisting in replacing the normal cone of a closed vector space by an arbitrary maximally monotone operator $C\colon\mathcal{H}\to 2^{\mathcal{H}}$, and consequently in equation 2.3 the orthogonal projector by the corresponding resolvent $J_{\gamma C}$. As shown in § 3 (in the augmented space and preconditioned setting), it is straightforward to establish that the preimage of $\mathrm{zer}(A+B+C)$ by $J_{\gamma C}$ is equal to the set of fixed points of

$$T_{\mathrm{FDR}} \stackrel{\mathrm{def}}{=} J_{\gamma A}\bigl(2J_{\gamma C}-\mathrm{Id}-\gamma B J_{\gamma C}\bigr) + \bigl(\mathrm{Id}-J_{\gamma C}\bigr)\ .$$

However, the convergence analysis is more delicate, since the latter operator cannot be factorized in a direct composition product analogous to the middle term in equalities equation 2.3, because $J_{\gamma C}$ is not a linear operator in general. This rather technical work has been carried out by Davis and Yin (2017), yielding that $T_{\mathrm{FDR}}$ is an $\alpha$-averaged operator whose constant depends on the step size and on the cocoercivity modulus of $B$ in the exact same way as shown by us for $T_{\mathrm{GFB}}$ (Raguet, 2014, III.2.3).

## 3 Generalized Forward-Douglas–Rachford with Preconditioning

We derive now the algorithmic scheme and the convergence analysis of the forward-Douglas–Rachford splitting algorithm, in the case of an arbitrary number of functionals in the split-ting problem P3. As we keep in mind convex optimization applications, we also consider the problem

$$\text{find } x \in \arg\min \sum_{i=1}^{n} g_i + f + h \, . \tag{P4}$$

Moreover, we substitute step sizes and splitting weigths by symmetric operators, as is now customary for preconditioning purpose. As we described along § 2.2, the forward-Douglas–Rachford is a generalized forward-backward splitting algorithm with a resolvent instead of a projection step. Actually, most of the following can be directly deduced from our work (Raguet and Landrieu, 2015) by replacing $\sum_{i=1}^{n} W_i z_i$ by $J_{\Gamma C}(\sum_{i=1}^{n} W_i z_i)$. We thus make extensive use of notations and results from the latter paper. However, for the sake of completeness, we first write in full the algorithmic scheme and our convergence results. Then, we explicit the necessary modifications for the convergence analysis.

### 3.1   Algorithmic Scheme and Convergence Results

Let us start by briefly recalling the main notations and hypothesis on operators and functionals considered. $\mathcal{H}$ is a real Hilbert space, and $\mathfrak{S}_+(\mathcal{H})$ and $\mathfrak{S}_{++}(\mathcal{H})$ are the set of bounded, linear selfadjoint operators over $\mathcal{H}$ that are *positive* and *strongly positive*, respectively. For $M \in \mathfrak{S}_{++}(\mathcal{H})$, we write $(x, y) \mapsto \langle x \mid y \rangle_M \overset{\text{def}}{=} \langle Mx \mid y \rangle$ and $x \mapsto \|x\|_M \overset{\text{def}}{=} \sqrt{\langle Mx \mid x \rangle}$ the inner product and norm induced by $M$ over $\mathcal{H}$, and by $\mathcal{H}_M$ the Hilbert space $\mathcal{H}$ endowed with this inner product. If $S \subseteq \mathcal{H}$ is a closed subspace, we denote $P_S$ (respectively $P_S^M$) the orthogonal projector over $S$ in $\mathcal{H}$ (respectively $\mathcal{H}_M$); if $x \in \mathcal{H}$, it is also convenient to denote $x^S \overset{\text{def}}{=} P_S x$. We also remind that the role of each subspace $\mathcal{H}_i$ in hypotheses H2 and h2 below is to better handle situations where the operator $A_i$ or the functional $g_i$ in the splitting depends only on a restricted subspace of $\mathcal{H}$.

(H1)   $B\colon \mathcal{H} \to \mathcal{H}$ has full domain, and $L \in \mathfrak{S}_{++}(\mathcal{H})$ is such that

$$\forall\, (x, y) \in \mathcal{H}^2, \quad \langle Bx - By \mid x - y \rangle \geq \|Bx - By\|_{L^{-1}}^2 \, .$$

(H2)   For each $i \in \{1, \ldots, n\}$, $A_i\colon \mathcal{H} \to 2^{\mathcal{H}}$ is maximally monotone in $\mathcal{H}$, and $\mathcal{H}_i \subseteq \mathcal{H}$ is a closed subspace such that $A_i = P_{\mathcal{H}_i} A_i P_{\mathcal{H}_i}$. Also, $C\colon \mathcal{H} \to 2^{\mathcal{H}}$ is maximally monotone.

(H3)   $\operatorname{zer}(\sum_{i=1}^{n} A_i + B + C) \neq \emptyset$.

We now formulate the analogous assumptions for the convex optimization problem P4.

(h1)   $f\colon \mathcal{H} \to \mathbb{R}$ is convex and everywhere differentiable such that its gradient in $\mathcal{H}_L$ is nonexpansive, where $L$ is defined in hypothesis H1.

(h2)   For each $i \in \{1, \ldots, n\}$, $g_i\colon \mathcal{H} \to\, ]-\infty, +\infty]$ is convex, proper and lower semicontinuous such that $g_i = g_i \circ P_{\mathcal{H}_i}$, where $\mathcal{H}_i$ is defined in hypothesis H2. Also, $h\colon \mathcal{H} \to\, ]-\infty, +\infty]$ is convex, proper and lower semicontinuous.

(h3)   Domain qualification and feasibility conditions:

     (i)   $0 \in \operatorname{sri}\{\operatorname{dom} h - \cap_{i=1}^{n} \operatorname{dom} g_i\}$ and $\forall\, i \in \{1, \ldots, n\}$, $0 \in \operatorname{sri}\{\operatorname{dom} g_i - \cap_{j=1}^{i-1} \operatorname{dom} g_j\}$,

     (ii)   $\arg\min \sum_{i=1}^{n} g_i + f + h \neq \emptyset$.

Then, we give the requirements on the preconditioners.

(C1) $\Gamma \in \mathfrak{S}_{++}(\mathcal{H})$ such that

    (i) $\|L^{1/2}\Gamma L^{1/2}\| < 2$, where $L$ is defined in hypothesis H1, and

    (ii) $\forall\, i \in \{1,\ldots,n\}$, $\Gamma(\mathcal{H}_i) \subseteq \mathcal{H}_i$, where $\mathcal{H}_i$ is defined in hypothesis H2.

(C2) For each $i \in \{1,\ldots,n\}$, $W_i \in \mathfrak{S}_+(\mathcal{H})$ is such that

    (i) $\ker W_i = \mathcal{H}_i^\perp$,

    (ii) $W_{i|\mathcal{H}_i} \in \mathfrak{S}_{++}(\mathcal{H}_i)$, and

    (iii) $\Gamma^{-1}W_i = W_i\Gamma^{-1}$.

Moreover,

    (iv) $\sum_{i=1}^n W_i = \mathrm{Id}_{\mathcal{H}}$.

We finally recall that if $g\colon \mathcal{H} \to\,]{-}\infty, +\infty]$ is a proper, convex, lower semicontinuous functional and $M \in \mathfrak{S}_{++}(\mathcal{H})$, we write $\mathrm{prox}_g^M$ for its proximity operator in $\mathcal{H}_M$, which is well-defined as $\mathrm{prox}_g^M\colon x \mapsto \arg\min_{\xi\in\mathcal{H}} \tfrac{1}{2}\langle x-\xi\,|\, M(x-\xi)\rangle + g(\xi)$. The algorithmic scheme is given in algorithm 1.

---

**Algorithm 1:** Preconditioned forward-Douglas–Rachford for monotone inclusion problem P3 under hypotheses H1–3; for convex optimization problem P4 under hypotheses h1–3, substitute $B$ with $\nabla f$, $J_{\Gamma C}$ with $\mathrm{prox}_h^{\Gamma^{-1}}$, and for all $i \in \{1,\ldots,n\}$, $J_{W_i^{-1}\Gamma A_i}$ with $\mathrm{prox}_{g_i|\mathcal{H}_i}^{\Gamma^{-1}W_i}$.

**Input:** $(z_i)_{1\le i\le n} \in \bigtimes_{i=1}^n \mathcal{H}_i$;   $\Gamma$, $(W_i)_{1\le i\le n}$ satisfying hypotheses C1–2;
    $\forall\, k \in \mathbb{N}$, $\rho_k \in\,]0, 2 - \tfrac{1}{2}\|L^{1/2}\Gamma L^{1/2}\|[$;
**Initialize:** $x \leftarrow J_{\Gamma C}(\sum_{i=1}^n W_i z_i)$;   $k \leftarrow 0$;
**repeat**
    $p \leftarrow 2x - \Gamma B x$;
    **for** $i \in \{1,\ldots,n\}$ **do**
        $z_i \leftarrow z_i + \rho_k\bigl(J_{W_i^{-1}\Gamma A_i}(p^{\mathcal{H}_i} - z_i) - x^{\mathcal{H}_i}\bigr)$;
    $x \leftarrow J_{\Gamma C}(\sum_{i=1}^n W_i z_i)$;
    $k \leftarrow k + 1$;
**until** *convergence*;
**return** $x$.

---

The following theorem ensures the convergence, and robustness to summable errors on the computations of each operator. For each iteration $k \in \mathbb{N}$, we denote by $b_k, c_k \in \mathcal{H}$ the errors when computing $B$ and $J_{\Gamma C}$, respectively, and for each $i \in \{1,\ldots,n\}$, by $a_{i,k} \in \mathcal{H}_i$ the error when computing $J_{W_i^{-1}\Gamma A_i}$.

**Theorem 3.1.** Set $(z_{i,0})_{1\le i\le n} \in \bigtimes_{i=1}^n \mathcal{H}_i$ and define $(x_k)_{k\in\mathbb{N}}$ the sequence in $\mathcal{H}$ together with $\bigl((z_{i,k})_{1\le i\le n}\bigr)_{k\in\mathbb{N}}$ the sequence in $\bigtimes_{i=1}^n \mathcal{H}_i$ such that for all $k \in \mathbb{N}$, $x_k = J_{\Gamma C}(\sum_{i=1}^n W_i z_{i,k}) + c_k$ and for all $i \in \{1,\ldots,n\}$,

$$z_{i,k+1} = z_{i,k} + \rho_k\bigl(J_{W_i^{-1}\Gamma A_i}\bigl(2x_k^{\mathcal{H}_i} - z_{i,k} - \Gamma(Bx_k + b_k)^{\mathcal{H}_i}\bigr) + a_{i,k} - x_k^{\mathcal{H}_i}\bigr), \tag{3.1}$$

where $b_k, c_k \in \mathcal{H}$, $a_{i,k} \in \mathcal{H}_i$, and $\rho_k \in\,]0, 2 - \tfrac{1}{2}\|L^{1/2}\Gamma L^{1/2}\|[$.
Under hypotheses C1–2 and H1–3, if

(i) $\sum_{k\in\mathbb{N}} \rho_k\left(2 - \frac{1}{2}\|L^{1/2}\Gamma L^{1/2}\| - \rho_k\right) = +\infty$, and

(ii) $\sum_{k\in\mathbb{N}} \rho_k \|b_k\| < +\infty$, $\sum_{k\in\mathbb{N}} \rho_k \|c_k\| < +\infty$ and $\forall\, i \in \{1,\ldots,n\}$, $\sum_{k\in\mathbb{N}} \rho_k \|a_{i,k}\| < +\infty$,

then the sequence $(x_k)_{k\in\mathbb{N}}$ defined above equation 3.1 converges weakly towards a solution of problem P3. If moreover

(iii) there exists $i \in \{1,\ldots,n\}$ such that $A_i$ is uniformly monotone,
or alternatively B or C is uniformly monotone,

then the convergence of $\left(x_k^{\mathcal{H}_i}\right)_{k\in\mathbb{N}}$, respectively $(x_k)_{k\in\mathbb{N}}$, is strong.

The following corollary specializes theorem 3.1 to the case of the convex optimization problem P4.

**Corollary 3.1.** *Let the sequence $(x_k)_{k\in\mathbb{N}}$ be defined by substituting, into equation 3.1, B with $\nabla f$, $J_{\Gamma C}$ with $\mathrm{prox}_h^{\Gamma^{-1}}$, and for all $i \in \{1,\ldots,n\}$, $J_{W_i^{-1}\Gamma A_i}$ with $\mathrm{prox}_{g_i|\mathcal{H}_i}^{\Gamma^{-1}W_i}$. If, in addition to hypotheses C1–2 and h1–3, assumptions (i) and (ii) of theorem 3.1 are satisfied, then the sequence $(x_k)_{k\in\mathbb{N}}$ converges weakly towards a minimizer of problem P4. If moreover there exists $i \in \{1,\ldots,n\}$ such that $g_i$ is uniformly convex, or alternatively f or h is uniformly convex, then $\left(x_k^{\mathcal{H}_i}\right)_{k\in\mathbb{N}}$, respectively $(x_k)_{k\in\mathbb{N}}$, converges strongly.*

### 3.2 Convergence Proof

It is convenient to work on the product space $\bigtimes_{i=1}^n \mathcal{H}_i$; we recall here the important definitions.

Let $\mathcal{H} \stackrel{\text{def}}{=} \bigtimes_{i=1}^n \mathcal{H}_i$, endowed with the inner product $\langle \boldsymbol{x}|\boldsymbol{y}\rangle \stackrel{\text{def}}{=} \sum_{i=1}^n \langle x_i \mid y_i \rangle$ and induced norm $\|\boldsymbol{x}\| \stackrel{\text{def}}{=} \sqrt{\langle \boldsymbol{x}|\boldsymbol{x}\rangle} = \sqrt{\sum_{i=1}^n \|x_i\|^2}$. The identity is $\mathbf{Id}$, and with the definitions from § 3.1, we define the following linear operators $\mathcal{H} \to \mathcal{H}$, $\boldsymbol{W}\colon \boldsymbol{x} \mapsto (W_i x_i)_{1\le i\le n}$ and $\boldsymbol{\Gamma}\colon \boldsymbol{x} \mapsto (\Gamma x_i)_{1\le i\le n}$. We also introduce $\boldsymbol{\Sigma}\colon \mathcal{H} \to \mathcal{H}\colon \boldsymbol{x} \mapsto \sum_{i=1}^n W_i x_i$ and $\boldsymbol{\Sigma}^*\colon \mathcal{H} \to \mathcal{H}\colon x \mapsto \left(x^{\mathcal{H}_i}\right)_{1\le i\le n}$, together with the set $\boldsymbol{S} \stackrel{\text{def}}{=} \mathrm{ran}\,\boldsymbol{\Sigma}^* = \{\boldsymbol{x} \in \mathcal{H} \mid \exists x \in \mathcal{H}, \forall\, i \in \{1,\ldots,n\}, x_i = x^{\mathcal{H}_i}\}$. Finally, we define $\boldsymbol{A}, \boldsymbol{C}\colon \mathcal{H} \to 2^{\mathcal{H}}$ and $\boldsymbol{B}\colon \mathcal{H} \to \mathcal{H}$ by

$$\boldsymbol{A} \stackrel{\text{def}}{=} \bigtimes_{i=1}^n A_{i|\mathcal{H}_i}, \quad \boldsymbol{B} \stackrel{\text{def}}{=} \boldsymbol{\Sigma}^* B \boldsymbol{\Sigma} \quad \text{and} \quad \boldsymbol{C} \stackrel{\text{def}}{=} \boldsymbol{W}\boldsymbol{\Sigma}^* C \boldsymbol{\Sigma} + N_{\boldsymbol{S}},$$

where the definition of the normal cone $N$ was given in the introduction § 2.2. We refer again the reader to our earlier paper (Raguet and Landrieu, 2015, §5.2) for the key properties of these objects. We recall in particular that $\boldsymbol{\Gamma}$ and $\boldsymbol{W}$ commute and $\boldsymbol{\Gamma}^{-1}\boldsymbol{W} \in \mathfrak{S}_{++}(\mathcal{H})$, and establish the following property.

**Proposition 3.1.** *With the above definitions, the following statements hold.*

(i) *If C is uniformly monotone in $\mathcal{H}$ with modulus $\phi$, then*

$$\forall\, (\boldsymbol{x},\boldsymbol{y}) \in \mathcal{H}^2,\ \forall\, (\boldsymbol{u},\boldsymbol{v}) \in \boldsymbol{W}^{-1}\boldsymbol{\Gamma}\boldsymbol{C}\boldsymbol{x} \times \boldsymbol{W}^{-1}\boldsymbol{\Gamma}\boldsymbol{C}\boldsymbol{y},\ \langle \boldsymbol{u} - \boldsymbol{v} \mid \boldsymbol{x} - \boldsymbol{y}\rangle_{\boldsymbol{\Gamma}^{-1}\boldsymbol{W}} \ge \phi(\|\boldsymbol{\Sigma}\boldsymbol{x} - \boldsymbol{\Sigma}\boldsymbol{y}\|).$$

(ii) $\boldsymbol{W}^{-1}\boldsymbol{\Gamma}\boldsymbol{C}$ *is maximally monotone in $\mathcal{H}_{\boldsymbol{\Gamma}^{-1}\boldsymbol{W}}$ and $J_{\boldsymbol{W}^{-1}\boldsymbol{\Gamma}\boldsymbol{C}} = \boldsymbol{\Sigma}^* J_{\Gamma C} \boldsymbol{\Sigma}$.*

*Proof.* Let $(x, y) \in \mathcal{H}^2$.

(i). Let $(u, v) \in W^{-1}\Gamma C x \times W^{-1}\Gamma C y$. By definition, there exists $(u, v) \in C\Sigma x \times C\Sigma y$ and $(u', v') \in N_S x \times N_S y$ such that $u = W^{-1}\Gamma(W\Sigma^* u + u')$ and $v = W^{-1}\Gamma(W\Sigma^* v + u')$. Thanks to the properties of $N_S$, $\langle u' - v' \mid x - y \rangle = 0$, hence we can develop $\langle u - v \mid x - y \rangle_{\Gamma^{-1}W} = \langle W\Sigma^* u - W\Sigma^* v \mid x - y \rangle = \sum_i \langle W_i(u - v) \mid x_i - y_i \rangle = \langle u - v \mid \sum_i W_i(x_i - y_i) \rangle = \langle u - v \mid \Sigma x - \Sigma y \rangle$. If $C$ is uniformly monotone in $\mathcal{H}$ with modulus $\phi$, by definition $\langle u - v \mid \Sigma x - \Sigma y \rangle \geq \phi(\|\Sigma x - \Sigma y\|)$.

(ii). If $C$ is only monotone, $\phi$ can be replaced by 0 in the above inequality, showing that $W^{-1}\Gamma C$ is monotone in $\mathcal{H}_{\Gamma^{-1}W}$. Suppose now that $y \in S$. According to lemma 5.3 (i) (ii) in the above paper (Raguet and Landrieu, 2015), $\Sigma\Sigma^* = \mathrm{Id}_\mathcal{H}$ and $\Sigma^*\Sigma_{|S} = \mathbf{Id}_S$. Thus using the definition of $J_{\Gamma C}$, we have the equivalences

$$y = \Sigma^* J_{\Gamma C}\Sigma x \iff \exists y \in \mathcal{H}, \; y = \Sigma^* y \text{ and } \Sigma x \in y + \Gamma C y,$$
$$\iff \Sigma^*\Sigma x \in y + \Sigma^*\Gamma C\Sigma y.$$

Now, according to proposition 5.1 (ii) in the same paper (Raguet and Landrieu, 2015), $\Sigma^*\Sigma = P_S^{\Gamma^{-1}W}$, so that $\Gamma^{-1}W(x - \Sigma^*\Sigma x) \in S^\perp$. With the definition of $N_S$ and $y \in S$, we get $x \in \Sigma^*\Sigma x + W^{-1}\Gamma N_S y$, and thus

$$y = \Sigma^* J_{\Gamma C}\Sigma x \iff x \in y + \Sigma^*\Gamma C\Sigma y + W^{-1}\Gamma N_S y,$$
$$\iff x \in y + W^{-1}\Gamma(W\Sigma^* C\Sigma + N_S)y,$$

where we substituted $\Sigma^*\Gamma = \Gamma\Sigma^* = W^{-1}W\Gamma\Sigma^* = W^{-1}\Gamma W\Sigma^*$ thanks to the properties of $\Sigma^*$, $\Gamma$ and $W$. This characterizes $J_{W^{-1}\Gamma C}$ with our definition of $C$. Altogether, $W^{-1}\Gamma C$ is monotone in $\mathcal{H}_{\Gamma^{-1}W}$, with a resolvent with full domain; it is thus maximally monotone by the theorem of Minty (1962). ∎

The generalized forward-Douglas–Rachford operator is based on the following fixed-point equation, which is then rewritten onto the product space $\mathcal{H}$.

**Proposition 3.2.** *Under hypotheses C1–2 and H2, $x \in \mathcal{H}$ is a solution of problem P3 if, and only if, there exists $(z_i)_{1 \leq i \leq n} \in \bigtimes_{i=1}^n \mathcal{H}_i$ such that $x = J_{\Gamma C}(\sum_{i=1}^n W_i z_i)$, and for all $i \in \{1, \ldots, n\}$,*

$$z_i = J_{W_i^{-1}\Gamma A_i}\left((2x - \Gamma B x)^{\mathcal{H}_i} - z_i\right) + \left(z_i - x^{\mathcal{H}_i}\right), \tag{3.2}$$

*that is to say if, and only if, there exists $z \in \mathcal{H}$ such that $x = J_{\Gamma C}\Sigma z$ and $z$ is a fixed point of the operator $T \colon \mathcal{H} \to \mathcal{H}$ defined by*

$$T \stackrel{\text{def}}{=} J_{W^{-1}\Gamma A}(2J_{W^{-1}\Gamma C} - \mathrm{Id} - \Gamma B J_{W^{-1}\Gamma C}) + (\mathrm{Id} - J_{W^{-1}\Gamma C}). \tag{3.3}$$

*Proof.* Let $x \in \mathcal{H}$. Following our former proof (Raguet and Landrieu, 2015, proposition 5.2), we have the equivalence

$$0 \in \sum_{i=1}^n A_i x + Bx + Cx \iff \exists (z_i)_{1 \leq i \leq n} \in \bigtimes_{i=1}^n \mathcal{H}_i, \begin{cases} \forall i, \; W_i(x - \Gamma Bx - z_i) \in \Gamma A_i x, \\ \text{and } \sum_{i=1}^n W_i z_i \in x + \Gamma Cx, \end{cases} \tag{3.4}$$

and it is easy to derive

$$W_i(x - \Gamma Bx - z_i) \in \Gamma A_i x \iff W_i\left((x - \Gamma Bx)^{\mathcal{H}_i} - z_i\right) \in \Gamma A_i x^{\mathcal{H}_i},$$
$$\iff (x - \Gamma Bx)^{\mathcal{H}_i} - z_i + x^{\mathcal{H}_i} \in W_i^{-1}\Gamma A_i x^{\mathcal{H}_i} + x^{\mathcal{H}_i},$$
$$\iff x^{\mathcal{H}_i} = J_{W_i^{-1}\Gamma A_i}\left((2x - \Gamma Bx)^{\mathcal{H}_i} - z_i\right),$$

leading to equation 3.2. The translation to the product space $\mathcal{H}$ is straightforward, using our description of the operators $J_{W^{-1}\Gamma A}$ (Raguet and Landrieu, 2015, proposition 5.1 (iii)) and $J_{W^{-1}\Gamma C}$ in proposition 3.1 (ii), and the definition of $B$. ∎

**Remark 3.1.** Following proposition 3.2, if $z \in \operatorname{fix} T$, with $x = J_{\Gamma C}\Sigma z$ the corresponding solution, the right-hand side of equivalence 3.4 yields on the product space $W(\Sigma^* x - \Gamma B \Sigma^* x - z) \in \Gamma A \Sigma^* x$ and $\Sigma^* \Sigma z \in \Sigma^* x + \Gamma \Sigma^* C x$, which can be combined as $\Sigma^* \Sigma z - z \in W^{-1}\Gamma A \Sigma^* x + \Gamma B \Sigma^* x + \Gamma \Sigma^* C x$. Now, as already pointed above in the proof of proposition 3.1 (ii), $\Gamma^{-1} W(z - \Sigma^* \Sigma z) \in \mathcal{S}^\perp$, and thus $\Sigma^* \Sigma z - z \in -W^{-1}\Gamma N_{\mathcal{S}}\Sigma^* x$, and one finally obtains $0 \in W^{-1}\Gamma A \Sigma^* x + \Gamma B \Sigma^* x + W^{-1}\Gamma C \Sigma^* x$. In other words, $\Sigma^* x = J_{W^{-1}\Gamma C} z \in \operatorname{zer}(W^{-1}\Gamma A + \Gamma B + W^{-1}\Gamma C)$.

Because of the nonlinearity of $J_{W^{-1}\Gamma C}$, analysis of the operator $T$ is more delicate than our generalized forward-backward operator of the form equation 2.3. As mentioned in the introduction, this can be done following a lemma due to Davis and Yin (2017). Convergence of a generalized forward-Douglas–Rachford ensues.

**Theorem 3.2.** Set $z_0 \in \mathcal{H}$ and define $(z_k)_{k \in \mathbb{N}}$ the sequence in $\mathcal{H}$ such that for all $k \in \mathbb{N}$,

$$z_{k+1} \stackrel{\text{def}}{=} z_k + \rho_k \Big( J_{W^{-1}\Gamma A}(2(J_{W^{-1}\Gamma C} z_k + c_k) - z_k - \Gamma B(J_{W^{-1}\Gamma C} z_k + c_k) + b_k) + a_k \\ - (J_{W^{-1}\Gamma C} z_k + c_k) \Big), \quad (3.5)$$

where $(a_k, b_k, c_k) \in \mathcal{H}^3$ and $\rho_k \in \,]0, 2 - \frac{1}{2}\|L^{1/2}\Gamma L^{1/2}\|[$.
Under hypotheses C1–2 and H1–3, if

(i) $\sum_{k \in \mathbb{N}} \rho_k \big(2 - \frac{1}{2}\|L^{1/2}\Gamma L^{1/2}\| - \rho_k\big) = +\infty$, and

(ii) $\sum_{k \in \mathbb{N}} \rho_k \|a_k\| < +\infty$, $\sum_{k \in \mathbb{N}} \rho_k \|b_k\| < +\infty$, and $\sum_{k \in \mathbb{N}} \rho_k \|c_k\| < +\infty$,

then there exists $z \in \operatorname{fix} T$ such that $x \stackrel{\text{def}}{=} J_{\Gamma C}\Sigma z$ is a solution of problem P3 and that

(a) $(z_k)_{k \in \mathbb{N}}$ converges weakly to $z$, and

(b) $(J_{\Gamma C}\Sigma z_k)_{k \in \mathbb{N}}$ converges weakly to $x$.

*If moreover*

(iii) *there exists $i \in \{1, \ldots, n\}$ such that $A_i$ is uniformly monotone,
or alternatively $B$ or $C$ is uniformly monotone,*

*then*

(c) $\left((J_{\Gamma C}\Sigma z_k)^{\mathcal{H}_i}\right)_{k \in \mathbb{N}}$ converges strongly to $x^{\mathcal{H}_i}$,
respectively $(J_{\Gamma C}\Sigma z_k)_{k \in \mathbb{N}}$ converges strongly to $x$.

*Proof.* (a). Proposition 3.1 (ii) shows that $W^{-1}\Gamma C$ is maximally monotone in $\mathcal{H}_{\Gamma^{-1}W}$, hence $J_{W^{-1}\Gamma C}$ is firmly nonexpansive with the characterization of Minty (1962). Similarly, we showed that $J_{W^{-1}\Gamma A}$ is also firmly nonexpansive, and that $\Gamma B$ is $\|L^{1/2}\Gamma L^{1/2}\|^{-1}$-cocoercive in $\mathcal{H}_{\Gamma^{-1}W}$ (Raguet and Landrieu, 2015, proposition 5.1 (iii) and (v), respectively). Altogether, proposition 3.1 of Davis and Yin (2017) provides that $T$ is $\alpha$-averaged, with $\alpha \stackrel{\text{def}}{=} (2 - \frac{1}{2}\|L^{1/2}\Gamma L^{1/2}\|)^{-1}$. Thus, there exists $R\colon \mathcal{H} \to \mathcal{H}$ nonexpansive in $\mathcal{H}_{\Gamma^{-1}W}$ such that $T = \alpha R + (1-\alpha)\mathbf{Id}$, and it is possible to rewrite iteration 3.5 as

$$z_{k+1} = z_k + \rho_k(Tz_k - z_k + d_k), \qquad (3.6)$$
$$= z_k + \alpha\rho_k(Rz_k + \alpha^{-1}d_k - z_k),$$

where for all $k \in \mathbb{N}$, $\|d_k\|_{\Gamma^{-1}W} \le \|a_k\|_{\Gamma^{-1}W} + \|b_k\|_{\Gamma^{-1}W} + (3 + \|L^{1/2}\Gamma L^{1/2}\|)\|c_k\|_{\Gamma^{-1}W}$; the term $\|a_k\|_{\Gamma^{-1}W}$ comes out only with triangle inequality, $\|b_k\|_{\Gamma^{-1}W}$ using also nonexpansivity of $J_{W^{-1}\Gamma A}$, and the coefficient in front of $\|c_k\|_{\Gamma^{-1}W}$ using also nonexpansivity of $J_{W^{-1}\Gamma C}$ and $\|L^{1/2}\Gamma L^{1/2}\|$-Lipschitzianity of $\Gamma B$. With (ii) and norms equivalence, we get $\sum_{k\in\mathbb{N}} \rho_k\|d_k\|_{\Gamma^{-1}W} < +\infty$. Moreover, for all $k \in \mathbb{N}$, $\alpha\rho_k < 1$ and thanks to (i), $\sum_{k\in\mathbb{N}}\alpha\rho_k(1-\alpha\rho_k) = \alpha^2\sum_{k\in\mathbb{N}}\rho_k(\alpha^{-1} - \rho_k) = +\infty$. Finally, proposition 3.2 and hypothesis H3 ensures fix $T \ne \emptyset$; but fix $T = $ fix $R$, and the results follows from Combettes (2004, lemma 5.1) together with proposition 3.2.

The rest of the proof is essentially a slightly simplified and improved version of the proof of theorem 2.1 of Davis and Yin (2017), identifying our operators $W^{-1}\Gamma A$, $\Gamma B$, $W^{-1}\Gamma C$ and parameters $(\rho_k)_{k\in\mathbb{N}}$ with their operators $A$, $C$, $B$, and parameters $(\lambda_k)_{k\in\mathbb{N}}$, respectively, and $\gamma \stackrel{\text{set}}{=} 1$. To lighten the notations, define $z_C \stackrel{\text{def}}{=} J_{W^{-1}\Gamma C}z$, and for all $k \in \mathbb{N}$,

$$z_{C,k} \stackrel{\text{def}}{=} J_{W^{-1}\Gamma C}z_k \qquad \text{and} \qquad z_{A,k} \stackrel{\text{def}}{=} J_{W^{-1}\Gamma A}(2z_{C,k} - z_k - \Gamma B z_{C,k})$$

so that $Tz_k = z_{A,k} + z_k - z_{C,k}$. The proof of the above lemma (Combettes, 2004, lemma 5.1) actually shows that $(Tz_k - z_k)_{k\in\mathbb{N}} = (z_{A,k} - z_{C,k})_{k\in\mathbb{N}}$ converges strongly to zero.

(b). By definition of the resolvent, we have the inclusions $z_k - z_{C,k} \in W^{-1}\Gamma C z_{C,k}$ and

$$u_k \stackrel{\text{def}}{=} (2z_{C,k} - z_k - \Gamma B z_{C,k}) - z_{A,k} \in W^{-1}\Gamma A z_{A,k}. \qquad (3.7)$$

Since $(z_k)_{k\in\mathbb{N}}$ is weakly convergent, it is bounded, and by nonexpansivity of $J_{W^{-1}\Gamma C}$, so is the sequence $(z_{C,k})_{k\in\mathbb{N}}$. It thus admits at least one weak cluster point. Let $y \in \mathcal{H}$ be such a weak cluster point, and $(k_j)_{j\in\mathbb{N}}$ be a sequence such that $z_{C,k_j} \xrightarrow[j\to+\infty]{} y$. Recalling $z_{A,k} - z_{C,k} \xrightarrow[k\to+\infty]{} 0$, we deduce $z_{A,k_j} \xrightarrow[j\to+\infty]{} y$. Similarly, the sequence $(\Gamma B z_{C,k_j})_{j\in\mathbb{N}}$ is bounded and admits at least one weak cluster point; up to extracting an other subsequence, we can assume that it is weakly convergent; in turn, $(u_{k_j})_{j\in\mathbb{N}}$ is also weakly convergent. Finally, $u_k + \Gamma B z_{C,k} + (z_k - z_{C,k}) = z_{C,k} - z_{A,k} \xrightarrow[k\to+\infty]{} 0$ and thus corollary 25.5 (iii) of Bauschke and Combettes (2011) applied to the sequences $u_{k_j} \in W^{-1}\Gamma A z_{A,k_j}$, $\Gamma B z_{C,k_j}$ and $(z_{k_j} - z_{C,k_j}) \in W^{-1}\Gamma C z_{C,k_j}$, $j \in \mathbb{N}$, shows in particular that the weak limit of the last one satisfies $(z - y) \in W^{-1}\Gamma C y$. But this characterizes $y = z_C$ uniquely, and thus $(z_{C,k})_{k\in\mathbb{N}}$ converges weakly to $z_C$. Recalling that $J_{W^{-1}\Gamma C} = \Sigma^* J_{\Gamma C}\Sigma$ and $\Sigma\Sigma^* = \text{Id}_{\mathcal{H}}$, applying the bounded linear operator $\Sigma$ to the above yields weak convergence of $(J_{\Gamma C}\Sigma z_k)_{k\in\mathbb{N}}$ to $x \stackrel{\text{def}}{=} J_{\Gamma C}\Sigma z$.

(c). Again by the definition of the resolvent, $z - z_C \in W^{-1}\Gamma C z_C$. In addition, remark 3.1 shows that $z_C \in \text{zer}(W^{-1}\Gamma A + \Gamma B + W^{-1}\Gamma C)$, and thus

$$u \stackrel{\text{def}}{=} z_C - z_k - \Gamma B z_C \in W^{-1}\Gamma A z_C. \qquad (3.8)$$

Altogether, monotonicity of $W^{-1}\Gamma A$, $\Gamma B$, and $W^{-1}\Gamma C$ ensures that

$$p_{A,k} \stackrel{\text{def}}{=} \langle z_{A,k} - z_C \mid u_k - u\rangle_{\Gamma^{-1}W} \ge 0, \qquad (3.9)$$
$$p_{B,k} \stackrel{\text{def}}{=} \langle z_{C,k} - z_C \mid \Gamma B z_{C,k} - \Gamma B z_C\rangle_{\Gamma^{-1}W} \ge 0, \text{ and} \qquad (3.10)$$
$$p_{C,k} \stackrel{\text{def}}{=} \langle z_{C,k} - z_C \mid (z_k - z_{C,k}) - (z - z_C)\rangle_{\Gamma^{-1}W} \ge 0, \qquad (3.11)$$

respectively. Now, develop

$$p_{A,k} + p_{B,k} + p_{C,k} = \langle z_{A,k} - z_{C,k} | u_k - u \rangle + \langle z_{C,k} - z_C | u_k + (z_k - z_{C,k}) + \Gamma B z_{C,k} \rangle,$$
$$= \langle z_{A,k} - z_{C,k} | u_k - u - (z_{C,k} - z_C) \rangle.$$

The right-hand operand of this last inner product stays bounded because it involves only Lipschitzian operators and $(z_k)_{k \in \mathbb{N}}$ is bounded. Since the left-hand operand converges strongly to zero as $k$ tends to infinity, so does the sum $p_{A,k} + p_{B,k} + p_{C,k}$; but since each term is nonnegative, they each converges in turn to zero.

By the definition of $A$, the inclusions 3.7 and 3.8 translate over the components of the product space $\mathcal{H}$ as, for all $i \in \{1, \ldots, n\}$, $(u)_i \in W_i^{-1} \Gamma A_i(z_C)_i$ and for all $k \in \mathbb{N}$, $(u_k)_i \in W_i^{-1} \Gamma A_i(z_{A,k})_i$. Suppose that there exists $i \in \{1, \ldots, n\}$ such that $A_i$ is uniformly monotone with modulus $\phi$. Then, the inequality 3.9 becomes $p_{A,k} = \sum_{j=1}^n \langle (z_{A,k})_j - (z_C)_j | (u_k)_j - (u)_j \rangle_{\Gamma^{-1}W_j} \geq \langle (z_{A,k})_i - (z_C)_i | (u_k)_i - (u)_i \rangle_{\Gamma^{-1}W_i}$, where, for all $j \in \{1, \ldots, n\} \smallsetminus \{i\}$, monotonicity of $W_j^{-1} \Gamma A_j$ in $\mathcal{H}_{\Gamma^{-1}W_j}$ has been used. This last inner product is $\langle (z_{A,k})_i - (z_C)_i | \Gamma^{-1} W_i(u_k)_i - \Gamma^{-1} W_i(u)_i \rangle$, hence the above inclusions and uniform monotonicity of $A_i$ finally yields $p_{A,k} \geq \phi(\|(z_{A,k})_i - (z_C)_i\|)$. But by definition, $\phi$ is nonnegative, nondecreasing and vanishes only at zero; with $p_{A,k} \xrightarrow[k \to +\infty]{} 0$, we deduce that $(z_{A,k})_i \xrightarrow[k \to +\infty]{} (z_C)_i = (\Sigma^* J_{\Gamma C} \Sigma z)^{\mathcal{H}_i} = x^{\mathcal{H}_i}$. But we have seen that $z_{C,k} - z_{A,k} \xrightarrow[k \to +\infty]{} 0$, so that in particular $(J_{\Gamma C} \Sigma z_k)^{\mathcal{H}_i} = (z_{C,k})_i \xrightarrow[k \to +\infty]{} x^{\mathcal{H}_i}$.

Alternatively, if $B$ is uniformly monotone with modulus $\phi$, our previous result (Raguet and Landrieu, 2015, proposition 5.1 (iv)) shows that inequality 3.10 becomes $p_{B,k} \geq \phi(\|\Sigma z_{C,k} - \Sigma z_C\|) = \phi(\|J_{\Gamma C} \Sigma z_k - x\|)$. The result follows again from the properties of $\phi$ and convergence of $(p_{B,k})_{k \in \mathbb{N}}$ to zero.

Finally, if $C$ is uniformly monotone with modulus $\phi$, then proposition 3.1 (i) shows that inequality 3.11 becomes $p_{C,k} \geq \phi(\|\Sigma z_{C,k} - \Sigma z_C\|)$, and the result follows as above. ∎

**Remark 3.2.** Strong convergence results with uniform monotonicity of an operator $A_i$ obviously improves on our previous results on the generalized forward-backward (Raguet et al., 2013; Raguet and Landrieu, 2015). Note also that, in contrast to theorem 2.1 of Davis and Yin (2017), our conclusions (b) and (c) do not require the parameters $(\rho_k)_{k \in \mathbb{N}}$ to be bounded away from 0 and $2 - \frac{1}{2}\|L^{1/2} \Gamma L^{1/2}\|$. Moreover, one can improve inequality 3.10 with cocoercivity of $\Gamma B$, as $p_{C,k} \geq \|L^{1/2} \Gamma L^{1/2}\|^{-1} \|\Gamma B z_{C,k} - \Gamma B z_C\|^2_{\Gamma^{-1}W}$. It is then easy to deduce that $\Gamma B z_{C,k} \xrightarrow[k \to +\infty]{} \Gamma B z_C$, and obtain strong convergence of $(z_{C,k})_{k \in \mathbb{N}}$ towards $z_C$ if $C$ is demiregular.

**Corollary 3.2.** *Theorem 3.1 and corollary 3.1 hold.*

*Proof.* Skipping some calculations, recursion equation 3.1 is a specific instance of equation 3.5, leading to theorem 3.1. Now, identifying the operators $B$, $C$, and each $A_i$ with, respectively, $\nabla f$, $\partial h$, and each $\partial g_i$, the derivation of hypotheses H1–3 from hypotheses h1–3 has been established in our previous work (Raguet and Landrieu, 2015, Corollary 5.1), proving corollary 3.1. ∎

## 4 Numerical Illustration

The aim of this section is twofold. First, it extends our previous numerical experiments (Raguet and Landrieu, 2015, § 4) in a different experimental setting, demonstrating the efficiency of the proposed algorithm on graph-structured optimization. Second, it compares the use of our generalized forward-backward and of the forward-Douglas–Rachford splitting algorithms, in the light of the differences described in § 1. To this regard, we regret that the authors of the so-called "three operators splitting scheme" (Davis and Yin, 2017) did not acknowledge the strong connection between their algorithm and ours, and hence failed to show the practical interest of their work. It would be interesting to run such comparison over the experiments they consider in their paper, but neither the source code nor the data are publicly available at the time of this writing. Therefore, we provide our own experimental settings, which are medium- and large-scale problems classically considered in signal processing or machine learning communities, and for which the use of the forward-Douglas–Rachford splitting algorithm seems especially relevant.

Each algorithm is carefully implemented in C++ with parallelization of most operators with OpenMP specifications, and run on a personal computer with eight cores at 2.40 GHz. The source code for PFDR is available at the author's GitHub repository,[1] source code of others are available upon request.

## 4.1 Brain Source Identification in Electroencephalography

In electroencephalography (EEG), the brain activity is recorded with high temporal resolution thanks to several electrodes located on the head of a subject. The goal of *brain source identification* in EEG is to retrieve the brain regions that were activated during the recording. The brain activity of the subject can be modeled as a set of *dipoles* situated in various brain regions, activated at different intensities. From the physics of the problem, it is possible to derive the linear operator linking the activity of the dipoles to the electrical potential recorded at the electrodes, called the *lead-field* operator.

### 4.1.1 Graph-Structured Regularized Inverse Problem

The number of electrodes (observations) is much smaller than the number of dipoles (brain regions). Moreover, some dipoles have highly correlated effects on the electrodes, while some others have weak influence. In consequence, inverting the lead-field operator is a difficult ill-posed problem. However, important priors can be taken into account, following for instance Gramfort et al. (2013); Becker et al. (2014). First, only a small subset of all brain regions are supposed to be active at the same time. Second, neighboring regions which are simultaneously active are likely to be active with similar intensities. Third, at some time points one may know that all active brain regions are activated with the same sign. These sparsity and positivity priors can be enforced by a model structured on a graph $G = (V, E)$, where each vertex of $V$ is a brain region and the edge set $E \subset V \times V$ contains each pair of spatially neighboring regions. The problem is to find an element of $\arg\min_{\mathbb{R}^V} F$, where for all $x \stackrel{\text{def}}{=} (x_v)_{v \in V} \in \mathbb{R}^V$,

$$F(x) \stackrel{\text{def}}{=} \tfrac{1}{2}\|y - \Phi x\|^2 + \sum_{(u,v)\in E} \lambda^{(\delta_1)}_{(u,v)} |x_u - x_v| + \sum_{v \in V} \lambda^{(\ell_1)}_v |x_v| + \iota_{\mathbb{R}_+^V}(x) \,, \qquad (4.1)$$

$y \in \mathbb{R}^N$ is the observation[2] over $N$ electrodes, and $\Phi \colon \mathbb{R}^V \to \mathbb{R}^N$ is the lead-field operator. The first term is an Euclidean norm ensuring coherence with the observation, the second term is a weighted *graph total variation* enforcing similarity between adjacent active brain regions, the third term is a weighted $\ell_1$-norm enforcing sparsity, and the forth term ensures positivity. The parameters $\lambda$ balance the influence of spatial regularity and sparsity, and should depend on the lead-field operator and on the noise over the observations.

---

[1] https://github.com/1a7r0ch3/CP_PFDR_graph_d1
[2] note that we use here both standard notations: when $N$ is an integer, $\mathbb{R}^N$ is the Cartesian product of $N$ copies of the set $\mathbb{R}$; when $E$ is a set, $\mathbb{R}^E$ is the set of all applications from $E$ to $\mathbb{R}$, isomorphic to $\mathbb{R}^{|E|}$.

Note that this is a spatial-only formulation, i.e. we consider the inverse problem over a single time point in the recording. Starting from a full spatiotemporal recording, we select the time point with highest activity (squared observed values summed over the electrodes). However, we use the full spatiotemporal information in order to estimate the noise level and the penalization parameters. Without delving into details, we use crude heuristics based on the *Stein's unbiased risk estimation*, as proposed by us earlier (Raguet, 2014, § V.3), and adapted to the current setting. In particular, one must take into account that the columns of the matrix representing $\Phi$ can be highly correlated, or have different norms; the latter requires parameters $\lambda$ varying along vertices and along edges.

### 4.1.2 Competing Algorithms

In spite of recent interest for the minimization of graph total variation, few algorithms allows to minimize efficiently equation 4.1 in its full generality. As far as we know, the most popular among practitioners is the primal-dual splitting algorithm of Chambolle and Pock (2011), so we include it in our tests for comparisons.

**Preconditioned generalized forward-backward splitting** (algorithm 1 with $h \stackrel{\text{set}}{=} 0$, PGFB). A natural way to catch functional equation 4.1 into a splitting of the form problem P2 is to set the first term of the functional as the smooth part, i.e. $f: x \mapsto \|y - \Phi x\|^2$, to split the graph total variation into $|E|$ different functionals $\{g_{(u,v)} \,|\, (u, v) \in E\}$, and consider the two remaining terms together in a last functional $g_{|E|+1}$.

We described earlier how to set the diagonal preconditioners and apply our algorithm in such case (Raguet and Landrieu, 2015, § 3 and § 4.2); consider also the two following differences with the former setting. First, the proximity operator of $g_{|E|+1}$ is now a soft-thresholding followed by a projector over the positive orthant. Second, the gradient of the smooth functional is now $\nabla f: x \mapsto \Phi^*(\Phi x - y)$, and its Hessian is $\nabla^2 f = \Phi^*\Phi$. For preconditioning purpose, we approximate the Hessian simply by its diagonal (also called *Jacobi approximation*). Moreover, we simply set the operator $L$ satisfying hypothesis h1 as $L \stackrel{\text{set}}{=} \ell \,\text{Id}$, where $\ell \stackrel{\text{set}}{=} \|\Phi^*\Phi\| = \|\Phi\|^2$ is estimated with power method. Also, we do not consider any *reconditioning* along iterations.

**Preconditioned forward-Douglas–Rachford splitting** (algorithm 1, PFDR). The form of the problem suggests to use the same splitting as above while particularizing the functional $h \stackrel{\text{set}}{=} g_{|E|+1}$ in problem P4. This is especially relevant regarding the first property discussed in the introduction: in contrast to PGFB above, at each iteration the iterate undergoes the proximity operator of $h$, ensuring sparsity and positivity.

**Preconditioned primal-dual algorithm of** Pock and Chambolle (2011, PPD). The functional is split as $F \stackrel{\text{def}}{=} g \circ \Lambda + h$, where we particularized in $h$ the last two terms of equation 4.1 just as

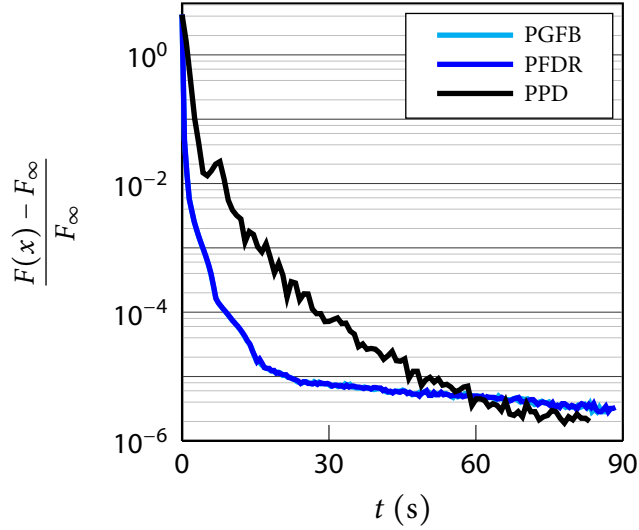

Figure 1: Brain source identification in EEG: optimization comparisons. PGFB and PFDR algorithms are so similar that one can hardly distinguish them.

for PFDR, and $g$ and $\Lambda$ are defined as

$$\Lambda: \begin{array}{ccc} \mathbb{R}^V & \longrightarrow & \mathbb{R}^N \times \mathbb{R}^E \\ x & \longmapsto & (v, \delta) \end{array}, \quad \text{with} \quad \left\{ \begin{array}{l} v \stackrel{\text{def}}{=} \Phi x, \text{ and} \\ \forall (u,v) \in E, \delta_{(u,v)} \stackrel{\text{def}}{=} \lambda_{(u,v)}^{(\delta_1)} (x_u - x_v), \end{array} \right.$$

$$\text{and} \quad g: \begin{array}{ccc} \mathbb{R}^N \times \mathbb{R}^E & \longrightarrow & \mathbb{R} \\ (v, \delta) & \longmapsto & \frac{1}{2} \|y - v\|^2 + \sum_{(u,v) \in E} |\delta_{(u,v)}| \end{array}.$$

Note that we use the preconditioned version of the algorithm; following the notations in the reference paper (our $g$, $h$ and $\Lambda$ being identified with their $F$, $G$ and $K$, respectively), the preconditioning matrices $T$ and $\Sigma$ are defined following lemma 2, equation (10), with the parameter $\alpha \stackrel{\text{set}}{=} 1$. We also set the relaxation parameter $\theta \stackrel{\text{set}}{=} 1$. Let us also mention that we tried other splitting and other parameters values, and that the above seem optimal for our purpose.

### 4.1.3 Decrease of the Functional and Brain Source Identification

The EEG data are provided by Ahmad Karfoul and Isabelle Merlet, LTSI, INSERM U1099. These are synthetic data; the graph $G$ and the lead-field operator $\Phi$ are computed according to a patient's model, then some brain sources are simulated, and finally the observations are constructed by application of the lead-field operator and of some noise. The model comprises 19 626 brain dipoles for only 91 electrodes. We select the time point of interest and the regularization parameters as described in § 4.1.1. Then, we run the algorithms with various stopping criteria and monitor the functional values against computational times. We also evaluate the quality of the solutions to assess the relevance of the algorithms for the problem at hand.

In order to set the best estimate of the solution, we perform $10^5$ iterations of each algorithm. After so many iterations all three give almost identical solutions and achieve almost equal functional value, set as the optimal value $F_\infty$. Then, difference between $F_\infty$ and functional values $F(x)$ along running time is given for each algorithm in figure 1. Results for PGFB and PFDR are so similar that one can hardly distinguish them. Comparing with PPD, one sees two regimes: at first, PGFB and PFDR achieves significantly lower functional values than PPD. Then for longer running times, the three behaves similarly. We believe that the initial speed-up of PGFB and PFDR is due to the use of the (preconditioned) gradient of the data-fidelity term, yielding reasonable reconstruction in a few iterations.

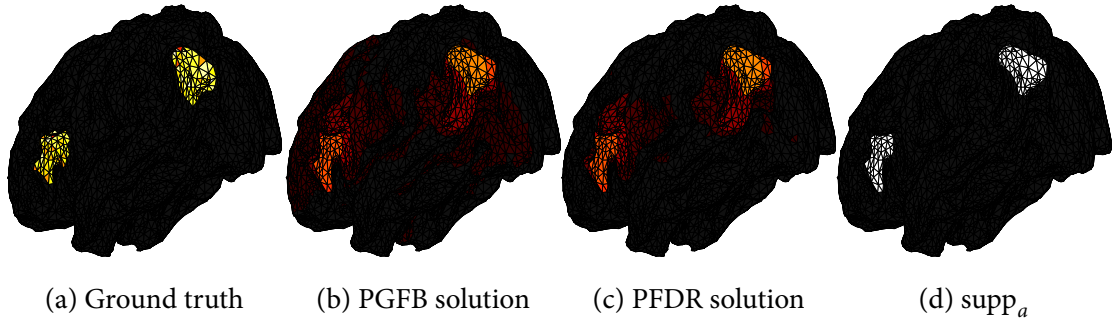

   (a) Ground truth      (b) PGFB solution     (c) PFDR solution       (d) $\operatorname{supp}_a$

Figure 2: Brain activity over the graph $G$: (a) ground truth, and retrieved with (b) PGFB and (c) PFDR (stopping criteria of $10^{-4}$ relative evolution, only 4 s running time). The same color map is used, from dark red (low activity) to bright yellow (high activity); dark gray indicates coefficients which are exactly zero. Note that lower amplitude of retrieved coefficients compared to ground truth is a typical effect of the convex regularizers we use. (d) Support recovered by discarding nonsignificant coefficients of solutions with 2-means clustering (see text); PGFB and PFDR yield exactly the same support.

Using synthetic data, the ground truth $\hat{x}$ is known (see figure 2(a)). Being interested in brain source identification, we consider the binary classification problem over the support $\operatorname{supp}(x) \stackrel{\text{def}}{=} \{v \in V \mid |x_v| > 0\}$, and evaluate the performance with the *Dice score* (also coined *F1 score*),

$$\operatorname{DS}(x;\hat{x}) \stackrel{\text{def}}{=} \frac{2|\operatorname{supp}(x) \cap \operatorname{supp}(\hat{x})|}{|\operatorname{supp}(x)| + |\operatorname{supp}(\hat{x})|}.$$

However, even the best minimizer of equation 4.1 (last line in table 1) gives poor prediction. This is mostly due to the fact that too many coefficients outside the real support have very low, but nonzero value. This can hardly be corrected by increasing the penalization parameters, because this results in bad reconstruction and eventually worst support recovery. In order to correct for this effect, we consider an approximated support $\operatorname{supp}_a(x) \stackrel{\text{def}}{=} \{v \in V \mid |x_v| > a\}$ where $a \in \mathbb{R}_+$ is determined by 2-means clustering of the absolute values of the solution $\{|x_v|\}_{v \in V}$.

As is customary in real experimental conditions, stopping criteria are prescribed as a minimum relative evolution of the iterate, $\|x_{\text{new}} - x_{\text{old}}\|/\|x_{\text{new}}\|$. A typical practitioner would set it according to a compromise between desired precision and available computational time. We propose several reasonable values, and for each one of them, we report the Dice score and the necessary running time of the algorithm in table 1.

As already explained, particularizing the $\ell_1$-norm as within PFDR and PPD schemes yields more zero coefficients than with PGFB, hence better raw predictions. Nevertheless, these are still far from satisfying (a Dice score around 0.3 is considered unreliable in any context). Now, discarding nonsignificant coefficients with a simple unidimensional 2-means clustering drastically improves the prediction (see also figures 2(b)–(d)). In that respect, PGFB behaves again exactly the same as PFDR, showing no practical advantage of the latter over the former. Finally, one can see that faster decrease of energy functionals of PGFB and PFDR in comparison to PPD translates somehow into prediction performance: the former algorithms yield good results in a few seconds, while the latter takes one order of magnitude longer for reaching the same level of prediction performance.

Table 1: Brain source identification: prediction performance and computing time comparisons.

|                        | PGFB |        |          | PFDR |        |          | PPD  |        |          |
|------------------------|------|--------|----------|------|--------|----------|------|--------|----------|
|                        | DS   | DS$_a$ | time (s) | DS   | DS$_a$ | time (s) | DS   | DS$_a$ | time (s) |
| rel. evol. $10^{-4}$   | 0.02 | 0.76   | 4        | 0.24 | 0.76   | 4        | 0.13 | 0.66   | 17       |
| rel. evol. $10^{-5}$   | 0.02 | 0.74   | 16       | 0.30 | 0.74   | 15       | 0.25 | 0.78   | 54       |
| rel. evol. $10^{-6}$   | 0.02 | 0.78   | 88       | 0.31 | 0.78   | 89       | 0.30 | 0.78   | 84       |
| $10^5$ iterations      | 0.04 | 0.78   | 258      | 0.32 | 0.78   | 245      | 0.32 | 0.78   | 251      |

## 4.2  Semantic Labeling of 3D Point Cloud

We consider now an application of our algorithm in the *remote sensing* field. Given a point cloud acquired from a LiDAR sensor, that is a list of points with their spatial coordinates in 3D, the goal is to assign a semantic label (or class, such as vegetation, building, car, etc…) to each point.

Recently, Guinard and Landrieu (2017, § 2.2) improve on the random forest classifier of Weinmann et al. (2015), using a combination of local features such as dimensionality and verticality (as described for instance by Demantké et al., 2011) and pre-estimated global features such as elevation and position with respect to the roads.

### 4.2.1  Convex Graph-Structured Spatial Regularization of Probabilistic Assignments

Such random forest allows to get relatively correct classifications with only a small number of training points, but this can be further improved by taking into account the spatial regularity of LiDAR point clouds: given the high spatial frequency of the sampling, neighboring points often belong to the same object. If $V$ denotes the set of points and $K$ the set of labels, the random forest classifier provides us with, for each point $v \in V$, a degree of confidence of the point belonging to each class under the form of a discrete probability distribution $(q_{v,k})_{k \in K} \in \triangle_K \overset{\text{def}}{=} \{ p \in \mathbb{R}^K \mid \sum_{k \in K} p_k = 1 \text{ and } \forall\, k \in K, p_k \geq 0 \}$. As advocated by us elsewhere (Landrieu et al., 2017), spatial regularity can again be enforced by minimizing a functional structured on a graph $G = (V, E)$, for instance defined for all $p \overset{\text{def}}{=} (p_{v,k})_{\substack{v \in V \\ k \in K}} \in \mathbb{R}^{V \times K}$ as

$$F(p) \overset{\text{def}}{=} \sum_{v \in V} \mathrm{KL}(\beta u + (1-\beta) q_v, \beta u + (1-\beta) p_v) + \sum_{(u,v) \in E} \lambda_{(u,v)} \sum_{k \in K} |p_{v,k} - p_{v,k}| + \iota_{\triangle_K^V}(p) \,, \quad (4.2)$$

where for all $r, s \in \triangle_K$, $\mathrm{KL}(r,s) \overset{\text{def}}{=} \sum_{k \in K} r_k \log(r_k / s_k)$ is the *Kullback–Leibler divergence*, $u \overset{\text{def}}{=} (1/|K|)_{k \in K} \in \triangle_K$ is the uniform discrete distribution, and $\beta \in\; ]0, 1[$ is a small smoothing parameter. The first term favors similarity with the original predictions; note that others are also considered, such as the opposite of the inner product between $p$ and $q$ which yields similar results in our experiments while favoring PPD algorithm in terms of computational time (data not shown). The second term is again a weighted graph total variation enforcing spatial regularity and the last term ensures that each labeling is a discrete probability distribution; The parameters $\lambda$ tune the influence of spatial regularity, and should depend on the confidence one has in the original predictions.

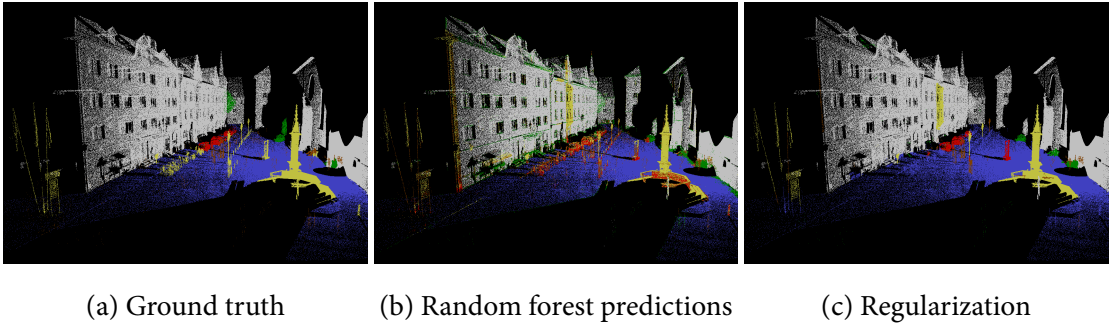

(a) Ground truth  (b) Random forest predictions  (c) Regularization

Figure 3: Semantic labeling of 3D point cloud. In our example, there are six different classes: road in blue, vegetation in green, façade in white, hardscape in yellow, cars in red, and scanning artifacts in brown. The improvement due to regularization is clear, but some flaws (like low hardscape items being treated as roads) suggest some refinements of the model, which could be implemented within the same framework but are beyond the scope of this article.

There are many ways of constructing the graph $G$ and selecting the penalization parameters; we briefly describe ours for the sake of completeness. We build the edge set simply by connecting each vertex to its ten nearest neighbors (with Euclidean distance). Then, we simply set the parameters $\left(\lambda_{(u,v)}\right)_{(u,v)\in E}$ constant along the edges. This constant is selected by a crude line search over a given range, where for each candidate penalization, an approximate solution $p$ of equation 4.2 is found, then each point $v \in V$ is assigned the label with its highest probability $\ell_v \in \arg\max_{k' \in K}\{p_{v,k'}\}$ (arbitrary chosen in the set if it is not a singleton), and this assignation is finally given a score computed on a training set $U \subset V$. The score we use is the *average of the F1 scores* across labels,

$$\overline{\mathrm{F}_1}(\ell; \hat{\ell}, U) \stackrel{\text{def}}{=} \frac{1}{|K|} \sum_{k \in K} \frac{2\left|\{v \in U \mid \ell_v = \hat{\ell}_v = k\}\right|}{\left|\{v \in U \mid \ell_v = k\}\right| + \left|\{v \in U \mid \hat{\ell}_v = k\}\right|} , \qquad (4.3)$$

where $\hat{\ell} \in K^V$ is the ground truth labelling. Let us finally describe an efficient selection of a small training set. The original predictions were obtained from a random forest classifier trained over a few tens of points within each class; the prediction given by $q$ over these points approaches perfect accuracy. Thus we complete this training set by adding, within each class, the same number of points, selected as the most uncertain; uncertainty being evaluated at each point $v \in V$ by the *entropy* of the prediction, defined as the quantity $\left(-\sum_{k \in K} q_{v,k} \log(q_{v,k})\right)$.

### 4.2.2 Competing Algorithms

We use the same proximal splitting algorithms as in § 4.1.2; some additional details are given below.

**Preconditioned generalized forward-backward splitting** (algorithm 1 with $h \stackrel{\text{set}}{=} 0$, PGFB). We proceed as previously, but now the smooth part in the splitting is the smoothed Kullback–Leibler divergence. Its gradient is $\nabla f \colon p \mapsto \left(-(1-\beta)\frac{\beta/|K|+(1-\beta)q_{v,k}}{\beta/|K|+(1-\beta)p_{v,k}}\right)_{\substack{v \in V \\ k \in K}}$ so that *hypothesis* h1 is satisfied

with the diagonal operator $L$ given by $\left((1-\beta)^2 \frac{\beta/|K|+(1-\beta)q_k}{(\beta/|K|)^2}\right)_{\substack{v\in V\\k\in K}}$. The proximity operator of $g_{|E|+1} \stackrel{\text{def}}{=} \iota_{\triangle_K^V}$ is an orthogonal projection over the simplex $\triangle_K$ independently for each vertex, which is relatively easy to compute in any diagonal metric. However, no heuristic is proposed to set the splitting weights $W_{|E|+1}$, which tune the importance of the simplex constraints along iterations. We simply set it to $0.2\,\text{Id}$ without more investigations, and set the other weights as already decribed, scaling them to sum up to $0.8\,\text{Id}$.

**Preconditioned forward-Douglas–Rachford splitting** (algorithm 1, PFDR). We particularize the functional $h \stackrel{\text{set}}{=} g_{|E|+1}$ in problem P4. This avoids completely to set splitting weight as above and ensures that each iterate belong to the product of simplices, in contrast to PGFB.

**Preconditioned primal-dual algorithm of** Pock and Chambolle (2011, PPD). The functional is split as $F \stackrel{\text{def}}{=} g \circ \Lambda + h$, where we particularized in $h$ the last term of equation 4.2 just as for PFDR, and $g$ and $\Lambda$ are defined as

$$\Lambda: \begin{array}{ccc} \mathbb{R}^{V\times K} & \longrightarrow & \mathbb{R}^{V\times K}\times\mathbb{R}^{E\times K} \\ p & \longmapsto & (p,\delta) \end{array} \quad \text{with } \forall\,(u,v)\in E,\, \delta_{(u,v),k} \stackrel{\text{def}}{=} \lambda_{(u,v)}(p_{u,k}-p_{v,k}),$$

$$g: \begin{array}{ccc} \mathbb{R}^{V\times K}\times\mathbb{R}^{E\times K} & \longrightarrow & \mathbb{R} \\ (p,\delta) & \longmapsto & \sum_{v\in V}\text{KL}(\beta u+(1-\beta)q_{v,k},\beta u+(1-\beta)p_{v,k}) + \sum_{(u,v)\in E}\sum_{k\in K}|\delta_{(u,v),k}|. \end{array}$$

The proximity operator of the smoothed Kullback–Leibler divergence can be easily computed in any diagonal metric, at the cost of finding a positive root of a second order polynomial for each coordinate.

### 4.2.3   Decrease of the Functional and Semantic Labeling

We use the dataset `domfountain1` from the database in http://www.semantic3d.net/, downsampled to 3 000 111 points for memory limitations. It is acquired with a fixed LiDAR, and is labeled with six different classes (see figure 3). The random forest are trained with only 25 points per classes, manually chosen; we then add 25 points per classes for selecting the penalization scaling as described in § 4.2.1. The smoothing parameter $\beta$ in the functional equation 4.2 is simply set to $0.1$ without further investigation.

We estimate the optimal functional value $F_\infty$ with $10^4$ iterations of PFDR; difference between $F_\infty$ and $F(p)$ along running time is given in figure 4. On this experiment, there are differences between PGFB and PFDR algorithms. First, PFDR enjoys a somewhat faster decrease of functional value along time. Second, missing functional values of PGFB indicates that, during the first iterations, some vertices of the iterate are so far away from the simplex that the Kullback–Leibler divergence have infinite value. However, they remain qualitatively similar, and both significantly faster than PPD which does not benefit from the smoothness of $f$.

Here, stopping criteria are prescribed as a minimum evolution of the iterate, $\|x_{\text{new}} - x_{\text{old}}\|_\infty$; we report the average F1 score described in equation 4.3 and the necessary running time of the algorithm in table 2. It is clear that PGFB and PFDR performs similarly, althought here PFDR gets faster to similar results. We believe that, as is the case for the brain source separation illustration, there exists weights $W_{|E|+1}$ which would make PGFB performs exactly as PFDR, but there is no easy way to set them optimally. Just as for the optimization consideration above, both are

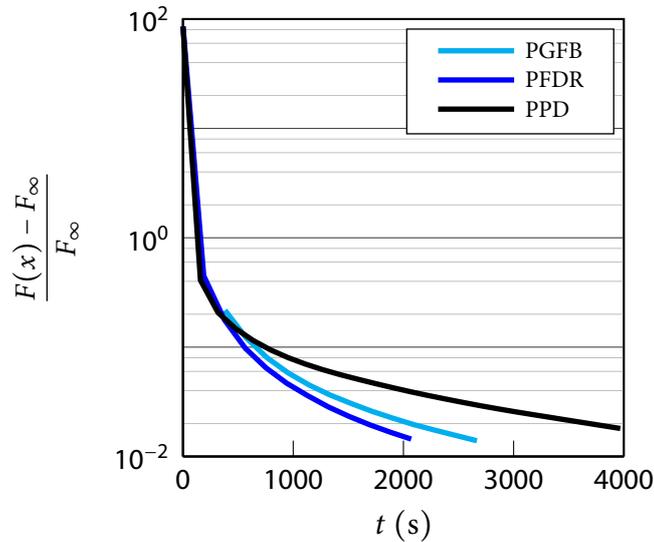

Figure 4: Semantic labeling of 3D point cloud: optimization comparisons.

Table 2: Semantic labeling: prediction performance and computing time comparisons. Average F1 score of random forest for comparison: 0.549.

|  | PGFB | | PFDR | | PPD | |
| --- | --- | --- | --- | --- | --- | --- |
|  | $\overline{F_1}$ | time (s) | $\overline{F_1}$ | time (s) | $\overline{F_1}$ | time (s) |
| rel. evol. $10^{-3}$ | 0.708 | 588 | 0.702 | 542 | 0.737 | 693 |
| rel. evol. $10^{-4}$ | 0.775 | 2706 | 0.776 | 2127 | 0.775 | 3968 |
| $10^4$ iterations | — | — | 0.776 | 9424 | — | — |

faster than PPD. Let us mention again for fairness that the latter is more adapted to the opposite of the inner product, in place of the smoothed Kullback–Leibler divergence, in the functional equation 4.2; which yields similar results for the problem considered here.

## 5 Conclusion

We write this note in the hope that the link between our generalized forward-backward and the forward-Douglas–Rachford gets clear, and that the latter gets this designation in its most general formulation. Additionally, we specify preconditioned and multiple splitting case, and give experimental evidence that its use in convex formulations of classical inverse or learning problems should be considered.

## References

H. H. Bauschke and P. L. Combettes. *Convex Analysis and Monotone Operator Theory in Hilbert Spaces.* Springer, New York, 2011.

H. Becker, L. Albera, P. Comon, R. Gribonval, and I. Merlet. Fast, variation-based methods for the analysis of extended brain sources. In *European Signal Processing Conference*, 2014.

L. M. Briceño-Arias. Forward-Douglas–Rachford splitting and forward-partial inverse method for solving monotone inclusions. *Optimization*, 64(5):1239–1261, 2015.

V. Cevher, C. B. Vu, and A. Yurtsever. Stochastic forward-Douglas–Rachford splitting for monotone inclusions. Technical report, EPFL, 2016.

A. Chambolle and T. Pock. A first-order primal-dual algorithm for convex problems with applications to imaging. *Journal of Mathematical Imaging and Vision*, 40(1):120–145, 2011.

P. L. Combettes. Solving monotone inclusions via compositions of nonexpansive averaged operators. *Optimization*, 53(5-6):475–504, 2004.

P. L. Combettes and J.-C. Pesquet. Primal-dual splitting algorithm for solving inclusions with mixtures of composite, Lipschitzian, and parallel-sum monotone operators. *Set-Valued and Variational Analysis*, 20(2):307–330, 2012.

L. Condat. A primal–dual splitting method for convex optimization involving lipschitzian, proximable and linear composite terms. *Journal of Optimization Theory and Applications*, 158 (2):460–479, 2013.

D. Davis and W. Yin. A three-operator splitting scheme and its optimization applications. *Set-Valued and Variational Analysis*, 25(4):829–858, 2017.

J. Demantké, C. Mallet, N. David, and B. Vallet. Dimensionality based scale selection in 3D LIDAR point clouds. *International Archives of the Photogrammetry, Remote Sensing and Spatial Information Sciences*, XXXVIII-5/W12:97–102, 2011.

J. Eckstein and D. P. Bertsekas. On the Douglas–Rachford splitting method and the proximal point algorithm for maximal monotone operators. *Mathematical Programming*, 55(3):293–318, 1992.

A. Gramfort, B. Thirion, and G. Varoquaux. Identifying predictive regions from fMRI with TV-$\ell_1$ prior. In *Pattern Recognition in Neuroimaging*. IEEE, 2013.

S. Guinard and L. Landrieu. Weakly supervised segmentation-aided classification of urban scenes from 3D LiDAR point clouds. *ISPRS Archives of the Photogrammetry, Remote Sensing and Spatial Information Sciences*, XLII-1/W1:151–157, 2017.

L. Landrieu, H. Raguet, B. Vallet, C. Mallet, and M. Weinmann. A structured regularization framework for spatially smoothing semantic labelings of 3D point clouds. *Journal of Photogrammetry and Remote Sensing*, 132:102–118, 2017.

P.-L. Lions and B. Mercier. Splitting algorithms for the sum of two nonlinear operators. *SIAM Journal on Numerical Analysis*, 16(6):964–979, 1979.

G. J. Minty. Monotone (nonlinear) operators in Hilbert space. *Duke Mathematical Journal*, 29: 341–346, 1962.

T. Pock and A. Chambolle. Diagonal preconditioning for first order primal-dual algorithms in convex optimization. In *IEEE International Conference on Computer Vision*, pages 1762–1769. IEEE, 2011.

H. Raguet, J. Fadili, and G. Peyré. A generalized forward-backward splitting. *SIAM Journal on Imaging Sciences*, 6(3):1199–1226, 2013.

H. Raguet. *A Signal Processing Approach to Voltage-Sensitive Dye Optical Imaging*. Ph.D. thesis, Université Paris-Dauphine, Paris, 2014.

H. Raguet and L. Landrieu. Preconditioning of a generalized forward-backward splitting and application to optimization on graphs. *SIAM Journal on Imaging Sciences*, 8(4):2706–2739, 2015.

J. E. Spingarn. Partial inverse of a monotone operator. *Applied Mathematics and Optimization*, 10(1):247–265, 1983.

B. C. Vũ. A splitting algorithm for dual monotone inclusions involving cocoercive operators. *Advances in Computational Mathematics*, 38(3):667–681, 2013.

M. Weinmann, B. Jutzi, S. Hinz, and C. Mallet. Semantic point cloud interpretation based on optimal neighborhoods, relevant features and efficient classifiers. *ISPRS Journal of Photogrammetry and Remote Sensing*, 105:286–304, 2015.

H. Zou, T. Hastie, and R. Tibshirani. On the "degrees of freedom" of the LASSO. *The Annals of Statistics*, 35(5):2173–2192, 2007.